\documentclass{amsart}
\usepackage{amssymb}

\newcommand{\REFER}{\ref}
\newcommand{\XREF}{\label}

\newcommand{\R}{{\mathbb R}}
\newcommand{\C}{{\mathbb C}}

\newcommand{\CC}{{\overline \C}}
\newcommand{\CR}{{\overline \R}}

\newcommand{\A}{{\mathcal A}}

\newcommand{\E}{{\mathcal E}}
\newcommand{\nbhd}{{neighborhood }}

\newcommand{\g}{\gamma}

\newtheorem{corollary}{Corollary}
\newtheorem{lemma}{Lemma}
\newtheorem{property}{Property}
\newtheorem{theorem}{Theorem}
\newtheorem{conjecture}{Conjecture}

\theoremstyle{definition}
\newtheorem{definition}{Definition}
\theoremstyle{remark}
\newtheorem{remark}{Remark}

\numberwithin{equation}{section}

\begin{document}

\title{Completely invariant sets of normality for rational semigroups}
\author{Rich Stankewitz}
\thanks{Research supported by a Department of Education 
GAANN fellowship and by the Research
Board of the University of Illinois at Urbana-Champaign.}
\address{Department of Mathematics, 
         University of Illinois, 
         Urbana, Illinois 61801}
\curraddr{Department of Mathematics, Texas A\&M University,
  College Station, TX 77843}
\email{richs@math.tamu.edu}

\begin{abstract}
Let $G$ be a semigroup of rational 
functions of degree at least two where the semigroup operation is
composition of functions.  We prove that the largest
open subset of the Riemann sphere on which the semigroup $G$ is normal
and is completely invariant under each element of $G$, 
can have only 0, 1, 2, or infinitely many components. 
\end{abstract}

\maketitle

\section{Introduction}

It is well known in iteration theory 
that the set of normality of a rational function can have
only 0, 1, 2, or infinitely many components (see~\cite{Be}, p.~94).  
In this paper we generalize this result by showing that 
the completely invariant set of normality of a
rational semigroup can have only 0, 1, 2, or infinitely many
components.  The proof not only
generalizes the iteration result, but it also provides an alternative proof
for it.

The results of this paper appear as part of the author's Ph.D. thesis,
written under the supervision of Aimo Hinkkanen.

\section{Definitions and basic facts}

In what follows all notions of convergence will be with respect to the 
spherical metric on the Riemann sphere $\CC.$

%\begin{definition}
%A family of meromorphic functions $\mathcal G$ is normal in a region
%$\Omega$ if every sequence of elements of $\mathcal G$ contains a
%subsequence which converges locally uniformly in $\Omega$ to a
%meromorphic function (possibly the constant $\infty$) defined on $\Omega$.
 
%\end{definition}

A rational semigroup $G$ is a semigroup 
of rational functions of degree greater than or equal to two
defined on the Riemann sphere $\CC$ with the semigroup operation being
functional composition.  When a semigroup $G$ is generated by the functions
$\{f_1, f_2, \dots, f_n, \dots\}$,
we write this as
\begin{equation}
G=\langle f_1, f_2,\dots, f_n, \ldots \rangle.\notag
\end{equation}

In~\cite{HM1}, p.~360 the definitions of the set of
normality, often called the Fatou set, 
and the Julia set of a rational semigroup are as follows:

\begin{definition} \label{N(G),J(G)}
For a rational semigroup $G$ we
define the set of normality of $G$, $N(G)$, by 
$$N(G)=\{z \in \CC:\text{there is a \nbhd of } z \text{ on
which } G \text{ is a normal family} \}$$
and define the Julia set of $G$, $J(G)$, by
$$J(G)=\CC \setminus N(G).$$
\end{definition}

Clearly from these definitions we see that $N(G)$ is an open set and
therefore its complement $J(G)$ is a compact set.
These definitions generalize the case of iteration of a single rational
function and we write $N(\langle h \rangle )=N_h$ and
$J(\langle h \rangle )=J_h$.  

Note that $J(G)$ contains the Julia set of each element of $G$.  In
fact, we have $J(G)=\overline{\cup_{f \in G} J_f}$ (see~\cite{HM1}, p.~365).

\begin{definition} If $h$ is a map of a set $Y$ into
itself, a subset $X$ of $Y$ is:
\begin{align} 
  i)& \,forward\,\,invariant \text{ under } h \text{ if } h(X) \subset X;\notag\\
 ii)& \,backward\,\,invariant \text{ under } h \text{ if }  
h^{-1}(X) \subset X;\notag\\
iii)& \,completely\,\,invariant \text{ under } h \text{ if } h(X) \subset
  X \text{ and } h^{-1}(X) \subset X.\notag
\end{align}

\end{definition}

It is well known that the set of normality of $h$ and the Julia set of
$h$ are
completely invariant under $h$ (see
~\cite{Be}, p.~54), in fact, 
\begin{equation}\label{cinv}
h(N_h)=N_h=h^{-1}(N_h) \text{ and }
h(J_h)=J_h=h^{-1}(J_h).
\end{equation}

Further, we have the following result.

\begin{property}\label{Jh2}
The set $J_h$ is the smallest closed completely invariant (under $h$)
set which contains three or more points (see ~\cite{Be}, p.~67).  
\end{property}
In fact, this may be chosen
as an alternate definition of $J_h$, equivalent to that 
given in Definition~\ref{N(G),J(G)}.

From Definition~\ref{N(G),J(G)}, we get that $N(G)$ is 
forward invariant under
each element of $G$ and $J(G)$ is backward invariant under each
element of $G$ (see ~\cite{HM1}, p.~360). 
The sets $N(G)$ and $J(G)$ are, however, not necessarily 
completely invariant under the elements of $G$.  
This is in contrast to the case of single function dynamics as noted
in (\ref{cinv}). 
The question then arises, what if we required 
the Julia set of the semigroup $G$ to be completely 
invariant under each element of $G$?  
We consider in this paper the consequences
of such an extension, given in the following
definition.

\begin{definition}\label{Edef}
 For a rational semigroup $G$ we define the completely invariant Julia
 set of $G$, $E=E(G)$, by 
$$E(G)=\bigcap\{S:S \text{ is closed, completely invariant under
each }g \in G, \#(S)\geq3 \}$$
where $\#(S)$ denotes the cardinality of $S$.
\end{definition}

We note that $E(G)$ exists, is closed, 
is completely invariant under each element
of $G$ and contains the Julia set of each
element of $G$ by Property~\REFER{Jh2}.

\begin{definition}\label{Wdef}
For a rational semigroup $G$ we define the completely invariant set of
normality of $G$, $W=W(G)$, to be the complement of $E(G)$, i.e., 
$$W(G)=\CC \setminus E(G).$$  
\end{definition}

Note that $W(G)$ is open and it is also completely invariant under each
element of $G$.
The main result of this paper is the following theorem.

\begin{theorem}\label{L<3}
For a rational semigroup $G$ the 
set $W(G)$ can have only 0, 1, 2, or infinitely many components.
\end{theorem}

\section{Proof of the main result}

\begin{property}\label{Eint}
If $E(G)$ has nonempty interior, then $E(G)=\CC.$
\end{property}

\begin{proof}
For a collection of sets $\A$, and a function $h$, we denote new
collections of sets by 
$h(\A)=\{h(A):A \in \A\}$ and $h^{-1}(\A)=\{h^{-1}(A):A \in \A\}$.

Choose $g \in G$.  Let us define the following countable collections of sets:
\begin{align}
&\E_0=\{ J_{g} \},\notag \\
&\E_1= \bigcup_{f \in G} f^{-1}(\E_0)  \cup \bigcup_{f \in G} f(\E_0),\notag\\
&\E_{n+1}=\bigcup_{f \in G} f^{-1}(\E_n) \cup  
\bigcup_{f \in G} f(\E_n),\notag\\
\text{ and }&\E = \bigcup_{n=0}^\infty \E_n.\notag
\end{align}

In the same manner as in the proof of Lemma 1
in~\cite{RS} we conclude $E(G)=\overline{\bigcup_{A \in \E} A}$.
Using this result we then finish the proof in the same manner as done  
in the proof of Lemma 2 in~\cite{RS}. 
\end{proof}

The remainder of this section will be devoted to the proof of
Theorem~\REFER{L<3}.

%Since $W$ is forward invariant under  $f$ and $g$ and $\CC
%\setminus W = E \supset J_f \cup J_g$ contains more than three points,
%we get $W \subset N(G) \subset
%N_f \cap N_g.$  From this we conclude that if $N_f$ or $N_g$ has 
%infinitely many components then either
%$W=\emptyset$ or $W$ has infinitely many components.

\begin{lemma}\label{compWtocompW}
If $W_0$ is a component of $W$, then $f(W_0)$ is
also a component of $W$ for any $f \in G$.
\end{lemma}

\begin{proof}
Let $W_1$ be the component of $W$ that contains $f(W_0)$.  We show
that $f(W_0)=W_1$.
Suppose to the contrary that $z \in W_1 \setminus f(W_0)$.  Since $f$
is continuous on the compact set $\overline{W_0}$ and an open map on
$W_0$, we have 
$\partial f(W_0) \subset f(\partial W_0) \subset f(E) \subset E$.  Let
$\gamma$ be a path in $W_1$ connecting $z$ to a point $w \in f(W_0)$.
Hence 
 $\gamma$ must cross $\partial f(W_0)
\subset E$.
This contradicts the fact that $\gamma \subset W_1$ and so we conclude
that $f(W_0)=W_1$.
\end{proof}

Since the remainder of this section will be devoted to the proof of
Theorem~\REFER{L<3}, we will assume that $W$
has $L$ components where $2 \leq L<+\infty.$  We remark here 
that the strategy will be to show that each of the $L$ components of
$W$ is simply connected and then the result will follow by an
application of the Riemann-Hurwitz relation.

\begin{definition} Let $W$ have components $W_j$ for $j=0,
\dots, L-1.$
\end{definition}

\XREF{fpermutes}
\begin{remark}
We see by Lemma~\REFER{compWtocompW} that each $f \in G$ (and hence each
$f^{-1}$ as well) permutes the
$W_j$ for $j=0, \dots, L-1$ since $f$ is a continuous map of $W$ onto
$W$.  
\end{remark}

We may assume that $\infty \in W_0$, else we may impose this condition
by conjugating each $f \in G$ by the same rotation of the sphere.

\begin{definition}\XREF{K's} For $j=1, \dots , L-1,$ we define
$$K_j=\{z \notin W_j: \text{there exists a simple closed curve }\gamma
\subset W_j \text{ such that } Ind_{\gamma}(z)=1\}$$
where the winding number is given by 
$Ind_{\gamma}(z)=(1/2\pi i)\int_{\gamma}1/(w-z)\,\, dw$.
If $z \in K_j$ and the simple closed curve $\gamma
\subset W_j$ is such that $Ind_{\gamma}(z)=1$, then we say that
$\gamma$ works for $z \in K_j$.
\end{definition}

In order to properly define $K_0$ we first need to move $W_0$ so that it
no longer contains $\infty.$
Let $\phi$ be a rotation of the sphere so that $\infty \in \phi(W_1)$ and
denote $\widetilde{W_j} = \phi(W_j)$ for $j=0, \dots, L-1$.

\begin{definition}\XREF{K0tilde} We define
$$\widetilde{K_0}=\{z \notin \widetilde{W_0}: \text{there exists a simple closed curve }\gamma
\subset \widetilde{W_0} \text{ such that } Ind_{\gamma}(z)=1\}$$
and
$$K_0= \phi^{-1}(\widetilde{K_0}).$$
If $z \in K_0$ and simple closed curve $\gamma
\subset \widetilde{W_0}$ is such that $Ind_{\gamma}(\phi(z))=1$, then we say that
the simple closed curve $\phi^{-1}(\gamma)$ works for $z \in K_0$.
\end{definition}

\begin{remark}\XREF{work}
Note that saying $\phi^{-1}(\gamma)$ works for $z \in K_0$ does
not necessarily imply that $Ind_{\phi^{-1}(\gamma)}(z)=1$, since it may
be the case that $Ind_{\gamma}(\phi(\infty))=1$ and hence
$Ind_{\phi^{-1}(\gamma)}(z)=0$ since $z$ lies in the unbounded
component of $\CC \setminus \phi^{-1}(\gamma).$
\end{remark}

\begin{definition} We define
$$K= \bigcup_{j=0}^{L-1} K_j.$$
\end{definition}

\begin{definition} We define
$$W_j'= W_j \cup K_j.$$
\end{definition}

\begin{lemma}\XREF{Wj'sc} 
For $j=0, \dots, L-1,$ the set $W_j'$ is open, connected and simply connected.
Thus each $K_j$ is the union of the ``holes'' in $W_j$.
\end{lemma}

\begin{proof}
Suppose that $1 \leq j \leq L-1$, so that 
$W_j$ is a bounded domain in the complex plane.
Define $A$ to be the unbounded component of $\CC \setminus W_j$.  
Hence $B=\CC \setminus A$ is open, connected and simply connected. 
  
Let $F$ be a bounded component of $\CC \setminus W_j$.  Since $A$ and
$F$ are each components of the closed set $\CC \setminus W_j$,
there exists a simple polygon $\g \subset W_j$ which separates $A$
from $F$ (see~\cite{Ne}, p.~134).  
Hence we see that $F \subset K_j$.  Since $F$ was an
arbitrary bounded component of $\CC \setminus W_j$, we conclude that
$K_j$ contains all the bounded components of $\CC \setminus W_j$,
i.e., the ``holes'' of $W_j$.  Hence $W_j' \supset B$.  Clearly $K_j$
cannot contain any points of $A$ since any simple closed path $\g \subset
W_j$ which would wind
around such a point would have to necessarily wind around every point
of $A$ (since $A$ is a component of the complement of $W_j$) including
$\infty$ which cannot happen.  Hence we conclude $W_j'=B$ and is
therefore open, connected and simply connected.

We show that $\phi(W_0')$ is open, connected and simply connected using
the same argument as above, and this implies that $W_0'$ is open,
connected and simply connected.
\end{proof}

\begin{definition} We define
$$W'= \bigcup_{j=0}^{L-1} W_j'.$$
\end{definition}

Note that we have $W'=W \cup K.$

\begin{lemma}\XREF{WrinWs} If for some distinct $r,s \in \{0, \dots, L-1\}$, 
we have $W_r' \cap W_s' \neq \emptyset$, then
 either $\overline{W_r'} \subset  W_s'$ or $\overline{W_s'} \subset  W_r'$.
In particular, if $W_r \cap W_s' \neq \emptyset$ for some distinct 
$r,s \in \{0, \dots, L-1\}$, 
then $\overline{W_r'} \subset  W_s'$.
\end{lemma}

\begin{proof}
Let $z \in W_r' \cap W_s'$.  Since $W_r \cap W_s = \emptyset$, we may
assume that $z \in K_s$, say.  Let $\g_s$ work for $z \in K_s$.
Let $I_{\g_s}$ be the component of $\CC \setminus \g_s$ which
contains $z$.  Note that $I_{\g_s} \setminus W_s= \{z:\g_s$ works for
$z\}$ whether or not $s=0$ (see Definitions~\REFER{K's}
and~\REFER{K0tilde} and Remark~\REFER{work}).  
Since $z \in W_r'$, we have two cases,
either $z \in K_r$ or $z\in W_r$.  

Suppose that $z \in K_r$ and let $\g_r$ work for $z \in K_r$.  As $\g_s
\cap \g_r =\emptyset$ (since $W_r \cap W_s = \emptyset$) we see that
either $\g_r \subset I_{\gamma_s}$ or $\g_s \subset I_{\gamma_r}$, 
where $I_{\g_r}$ is
the component of $\CC \setminus \g_r$ which
contains $z$.  By switching the
roles of $r$ and $s$, if necessary, we assume $\g_r \subset I_{\gamma_s}$ and
we note that this can be done since $z \in K_r \cap K_s$.  
In particular,
$W_r \cap I_{\gamma_s} \neq \emptyset$.  

If $z \in W_r$, then we still get 
$W_r \cap I_{\gamma_s} \neq \emptyset$ since $z \in I_{\gamma_s}$.  

Since $W_r \cap I_{\gamma_s} \neq \emptyset$, $W_r \cap W_s = \emptyset$,
$\overline{W_r}$ is connected, and
$\g_s \subset W_s$, we conclude that $\overline{W_r} \subset
I_{\gamma_s}.$  Hence 
$\overline{W_r} \subset W_s'$ since $\g_s$ then works for every
$z \in \overline{W_r}$.  Since $W_s'$ is simply connected we
see that $\overline{W_r'} \subset  W_s'$.
\end{proof}

\begin{lemma}\XREF{W0'} 
The boundary of $W_0'$ is a nondegenerate continuum and as such
contains more than three points.
\end{lemma}

\begin{proof}
We will first show that $W_0' \cap W_1' = \emptyset$.  The set $W_1'$
cannot contain $W_0'$ as $\infty \in W_0'$ and $W_1'$ is a bounded
subset of $\C$ (since $W_1$ is a bounded subset of $\C$).
The same argument also shows that $\phi(W_0')$ cannot contain
$\phi(W_1')$ where $\phi$ is as in Definition~\REFER{K0tilde}, 
and so we conclude that $W_0'$
cannot contain $W_1'$.
By Lemma~\REFER{WrinWs} we conclude that $W_0' \cap W_1' = \emptyset$.

Since $W_0'$ is simply connected, $\partial W_0'$ contains a
nondegenerate continuum unless 
$\partial W_0'$ consists of just a single point.  If 
$\partial W_0'$ consists of just a single point, then $W_0' \cup
\partial W_0' = \CC$, but this contradicts the fact that 
$W_0' \cap W_1' = \emptyset$.
\end{proof}

\begin{lemma}\XREF{Wjbdry}
For each $j=0,\dots, L-1,$ we have $J_f \subset \partial W_j$ for each
$f \in G$.  Since $J(G)=\overline{\cup_{f \in G} J_f}$, we have
$J(G) \subset \partial W_j$ for each $j=0,\dots, L-1$.
\end{lemma}

\begin{proof}
Since $f$ permutes the $W_j$ by Remark~\REFER{fpermutes}, we may select 
a positive integer $n$ so that $f^n(W_j)=W_j=f^{-n}(W_j)$ 
for each $j=0,\dots,
L-1.$   Then we have $\overline{\bigcup_{k=1}^\infty f^{-kn}(W_j)}
\supset J_{f^n}=J_f$ (see~\cite{Be}, p.71 and p.51).  But since 
$\overline{\bigcup_{k=1}^\infty f^{-kn}(W_j)} = \overline{W_j}$ we see
that $\partial W_j \supset J_f$, since $W_j \cap J_f = \emptyset$.
\end{proof}

\begin{lemma}\XREF{WrnotinWs} 
We have $W_r \nsubseteq  W_s'$ for distinct 
$r,s \in \{0, \dots, L-1\}$, 
and therefore, by Lemma~\REFER{WrinWs}, 
the $W_j'$ are disjoint for $j=0, \dots, L-1$.  
\end{lemma}

\begin{proof}
If $L=2$, then the proof of Lemma~\REFER{W0'} shows that $W_0' \cap W_1'
= \emptyset$.

We assume now that $L \geq 3$.  
We will first show that no bounded $W_s'$ can contain any $W_r$ with
$r \neq s.$  
Suppose that this does occur.  Then there exists a simple closed 
curve $\g_s \subset W_s$
such that $\overline{W_r} \subset I_{\g_s}$ where $I_{\g_s}$ 
is the component of
$\CC\setminus \g_s$ which contains the points $z$ such that
$Ind_{\g_s}(z) = 1.$  Hence, by Lemma~\REFER{Wjbdry},
$J(G) \subset \partial W_r\subset 
\overline{W_r} \subset I_{\g_s}$.  
But since $\overline{W_0} \subset \CC \setminus
I_{\g_s}$ we see that $J(G) \subset \partial W_0\subset 
\overline{W_0} \subset \CC \setminus I_{\g_s}$.  This contradiction implies 
no bounded $W_s'$ can contain any $W_r.$

We see that $W_0'$ cannot contain any $W_r$ with $r \geq 1$
by the following similar argument.  If $\overline{W_r}\subset W_0'$,
then there exists a simple closed curve $\g \subset \widetilde{W_0}$ such
that $Ind_{\g}(z)=1$ for every $z \in \overline{\widetilde{W_r}}$.  Let
$I_{\g}$ be the component of $\CC \setminus \g$ which contains
$\overline{\widetilde{W_r}}$.  So $\phi(J(G)) \subset
\phi(\partial W_r)=\partial \phi(W_r)=\partial \widetilde{W_r}\subset
I_{\g}.$  Since $\overline{\widetilde{W_1}} \subset \CC \setminus I_{\g}$
(recall $\infty \in \widetilde{W_1}$), we see that $\phi(J(G)) \subset
\phi(\partial W_1)=\partial \phi(W_1)=\partial \widetilde{W_1} \subset \CC
\setminus I_{\g}.$  This contradiction implies $W_0'$ cannot contain
any $W_r$ with $r \geq 1$.
\end{proof}

\begin{corollary}\XREF{Knoint} The 
set $K$ has no interior and therefore
each $K_j \subset \partial W_j$.
\end{corollary}

\begin{proof}
By Lemma~\REFER{WrnotinWs} we see that each $K_j \subset E$ and hence
$K \subset E$.  The Corollary then follows from Property~\REFER{Eint}.
\end{proof}

\begin{corollary}\XREF{bdryWj} We have 
$\partial W_j = K_j \cup \partial W_j'$.
\end{corollary}

\begin{proof}
By Corollary~\REFER{Knoint} we get $K_j \cup \partial W_j'
\subset \partial W_j.$  We also have $\partial W_j=\overline{W_j}
\setminus W_j \subset \overline{W_j'}
\setminus W_j =(W_j' \cup \partial W_j') \setminus W_j=
(W_j\cup K_j \cup \partial W_j') \setminus W_j= K_j \cup \partial
W_j'.$  
\end{proof}

\begin{lemma} \XREF{Kfinv}
We have $f(K) \subset K$ for all $f \in G$.
\end{lemma}

\begin{proof}
Let $z \in K_j$ be such that $\g \subset W_j$ works for $z$.

Suppose that $W_l=f(W_j) \neq W_0$.  So $W_j'$ contains no poles of $f$,
else such a pole would be in $W_j$ (by the complete invariance of $W$
under $f$ since $\infty \in W_0 \subset W$ and Lemma~\REFER{WrnotinWs}) and
hence $f(W_j) = W_0$.  By the argument principle, $f(\g)\subset W_l$ winds
around $f(z)$, thus $f(z) \in K_l$ as $f(z) \notin W_l$ by the
complete invariance of $W$ under the map $f$. Note that $f(\g)$ might
not work for $f(z) \in K_l$ since it might not be simple, but
$f(z) \in K_l$ since it cannot be in the unbounded component of $\CC
\setminus W_l$ and have a curve in $W_l$, namely $f(\g)$, wind around
it. 

Now suppose that $f(W_j) = W_0$.  So $(\phi \circ f)(W_j) = \widetilde{W_0}$
is bounded and $W_j'$ contains no poles of $\phi \circ f$ (else
$f(W_j) = W_1$).  So $(\phi \circ f)(\g)$ winds around $(\phi \circ f)(z)$
and hence $(\phi \circ f)(z) \in \widetilde{K_0}$, i.e., $f(z) \in K_0$.

So $f(K_j) \subset K$ and hence we conclude $f(K) \subset K.$  
\end{proof}

\begin{lemma} \XREF{W'inNfandNg}
We have for all $f \in G$, $f(W')\cap \partial W_0'=\emptyset$.  
Also $W' \subset N(G)$ and in particular $K \cap J(G) = \emptyset.$
\end{lemma}

\begin{proof}
We have $f(W')=f(W \cup K)=f(W) \cup f(K) \subset W \cup K=
W'$.
Since $W' \cap \partial W_0' = \emptyset$ (since $W'$ is open),
Lemma~\REFER{W0'} and Montel's Theorem finish the proof.
\end{proof}

\begin{corollary}\XREF{bdryWj'}
We have $J(G) \subset \partial W_j'$ 
for each  $j=0, \dots, L-1$.  
\end{corollary}

\begin{proof}
This follows immediately from Lemma~\REFER{Wjbdry},
Corollary~\REFER{bdryWj} and Lemma~\REFER{W'inNfandNg}.
\end{proof}

\begin{remark}\XREF{Wada}
It is of interest to note that for any
positive integer $n$ there exist disjoint simply connected domains $D_1,
\dots, D_n$ in $\CC$ with $\partial D_1= \partial D_2 = \dots =
\partial D_n$ (see~\cite{HY}, p.~143).
Thus Corollary~\REFER{bdryWj'} does not imply that $L<3$ from a purely
topological perspective.
\end{remark}

\begin{lemma} \XREF{Kbinv}
We have $f^{-1}(K) \subset K$ for all $f \in G$.  Hence by
Lemma~\REFER{Kfinv}, $K$ is completely invariant under each
$f \in G$.
\end{lemma}

\begin{proof}
Let $z \in K_j \subset \partial W_j$ and say $f(w)=z$.  Define $W_k =
f^{-1}(W_j)$ by Remark~\REFER{fpermutes}.  We obtain
sequences $z_n \in W_j$ such that $z_n \to z$, and 
$w_n \in W_k$ such that $w_n \to w$ and $f(w_n)=z_n$.  Hence we see that $w \in
\partial W_k$, else $w \in W_k$ and $z=f(w) \in W_j$.  If $w \notin
K_k$, then $w \in \partial W_k'$ by Corollary~\REFER{bdryWj}.  Let
$\Gamma$ be the component of $\partial W_j$ that contains $f(\partial
W_k')$.  Since $z \in \Gamma$, the set $\Gamma$ must be one of the 
components of $K_j.$  By Corollary~\REFER{bdryWj'} we see that
there exists a $\zeta \in \partial W_k' \cap J_f$.  Hence $f(\zeta) \in K_j
\cap J_f$ which is a contradiction since we know by
Lemma~\REFER{W'inNfandNg} that $K$ is disjoint from $J(G) \supset J_f$.  
This contradiction implies $w \in K_k$ and hence $f^{-1}(K) \subset K$.
\end{proof}

\begin{lemma} \XREF{sc}
If $W$ has $L$ components where $2 \leq L<+\infty$, 
then each is simply connected.
\end{lemma}

\begin{proof}
Since $K$ and $W$ are each completely invariant under each $f \in G$, so is
$W'=W \cup K$.  By Lemma~\REFER{W'inNfandNg} we see that $\CC
\setminus W'$ is completely invariant under each $f \in G$, 
closed, and contains $J(G)$.  
Hence $E \subset \CC \setminus W'$.  This implies that $W=W'$
and hence each component of $W$ is then simply connected.
\end{proof}

We are now able to present the proof of Theorem~\REFER{L<3}.

\begin{proof}[Proof of Theorem~\REFER{L<3}]
If $W$ has $L$ components where $2 \leq L<+\infty$, 
then each is simply connected by Lemma~\REFER{sc}.
Select a map $f \in G$.  
Letting $n \geq 1$ be selected so that each of the components $W_j$ of $W$
is completely invariant under $f^n$, we get by the Riemann-Hurwitz
relation (see~\cite{S}, p.~7)
 
$$\delta_{f^n}(W_j)=\deg(f^n)-1$$
where we write $\delta_{g}(B)=\sum_{z \in B} [v_g(z)-1]$
and $v_g(z)$ is the valency of the map $g$ at the point $z$.

Hence we obtain
$$L(\deg(f^n)-1) = \sum^{L-1}_{j=0} \delta_{f^n}(W_j) 
\leq \delta_{f^n}(\CC) =2(\deg (f^n) -1)$$
and so $L \leq 2.$  The last equality follows from Theorem 2.7.1
in~\cite{Be}.
\end{proof}

\XREF{L=2}
\begin{remark} Note that if $L=2$, then each component of
$W$ is necessarily simply connected.
\end{remark}

\section{Conclusions}

We know from iteration theory that each of the four possibilities $(0,
1, 2, \infty)$ for the number of components of the set of normality can
be achieved.  So by constructing semigroups $G$ such that all the
elements have the same Julia set we know that the only four
possibilities for the number of components of the completely invariant
set of normality of $G$ can also be achieved.  However,  it does not seem
possible that all four possibilities can be achieved
if we restrict ourselves to the cases where 
two elements of the semigroup $G$ have nonequal Julia sets.  For
example, if $G$ contains two polynomials with nonequal Julia sets
then the completely invariant set of normality is necessarily empty
(see~\cite{RS}, Theorem 1).  

We do have the following examples however.

\vskip.1truein
{\bf{Example 1}}
Let $f(z)=2z - z^{-1}$ and $g(z)=(z^2 -1)/2z.$  In this
case $J_f$ is a Cantor subset of the interval $[-1,1]$ and
$J_g=\CR$, the extended real line.  It is shown
in~\cite{RS2} that $J(\langle f,g \rangle)=\CR=E(\langle f,g \rangle)$.

\vskip.1truein
{\bf{Example 2}}
Let $f(z)=2z - z^{-1}$ and $g(z)=f(z-1)+1.$  In this
case $J_f$ is a Cantor subset of the interval $[-1,1]$ and
$J_g=J_f+1$.  It is shown
in~\cite{RS2} that $J(\langle f,g \rangle)=[-1,2]$ and 
$E(\langle f,g \rangle)=\CR$.

\vskip.1truein
So we see that it is possible for a completely invariant set of
normality of a semigroup $G$ which contains two elements with nonequal
Julia sets, to have 0 or exactly 2 components.  
We feel that the interplay between functions with nonequal
Julia sets and the fact that if $E(G)$ has interior then
$E(G)=\CC$ demands that only under special circumstances can we have
$W(G)$ be nonempty, when two elements of the semigroup $G$ have
nonequal Julia sets.  

We state the following conjectures which are due to Aimo Hinkkanen and
Gaven Martin.

\begin{conjecture}\label{conj1}
If $G$ is a rational semigroup which contains two maps $f$ and $g$ such that
$J_f \neq J_g$ and $E(G) \neq \CC$, then $W(G)$ has exactly two components,
each of which is simply connected, and $E(G)$ is equal to the boundary of
each of these components.
\end{conjecture}

\begin{conjecture}\label{conj2}
If $G$ is a rational semigroup which contains two maps $f$ and $g$ such that
$J_f \neq J_g$ and $E(G) \neq \CC$, then $E(G)$ is a simple closed curve in
$\CC$.
\end{conjecture}

Of course Conjecture~\REFER{conj1} would follow from Conjecture~\REFER{conj2}.

We finish by including some comments on the number of components of
the set of normality $N(G)$ of a rational semigroup $G$.
It is not known if the set $N(G)$ must 
have only 0, 1, 2, or infinitely many components
when $G$ is a finitely generated rational semigroup.
However, for each positive integer $n$, 
an example of an infinitely generated polynomial
semigroup $G$ can be constructed with the property that $N(G)$ has
exactly $n$
components.  These examples were constructed by David Boyd in~\cite{Bo}.

\bibliographystyle{amsplain}
\bibliography{numcomp}

\enddocument